\documentclass[11pt, a4paper]{article}
\usepackage{setspace}
\usepackage{t1enc}
\usepackage[latin1]{inputenc}
\usepackage[english]{babel}
\usepackage{amssymb}
\usepackage{amsmath,amsthm}
\usepackage{amsfonts}
\usepackage{latexsym}
\usepackage[dvips]{graphicx}
\usepackage{float}
\usepackage[natural]{xcolor}
\usepackage{algorithm}
\usepackage{algorithmic}
\usepackage{enumerate}
\usepackage{multirow}
\usepackage{xcolor}
\usepackage[colorlinks,linkcolor=blue]{hyperref}
\usepackage{lineno}
\usepackage[a4paper, margin=1in]{geometry}
\usepackage{bm}
\usepackage{comment}
\usepackage{tikz}
\usepackage{pgfplots}
\usetikzlibrary{calc}
\usepackage{cleveref}
\usepackage{mathtools}
\usepackage{cite}
\usepackage{listings}
\title{Supersaturation of odd linear cycles}

\author
{Lirong Deng\footnote{School of Mathematics and Statistics, Beijing Institute of Technology. Email: {\tt lirong.deng@bit.edu.cn}.}
\and  Jie Han\footnote{School of Mathematics and Statistics, Beijing Institute of Technology. Email: {\tt han.jie@bit.edu.cn}. Research partially supported by National Natural Science Foundation of China (12371341).} \and Jiaxi Nie\footnote{School of Mathematics, Georgia Institute of Technology {\tt jnie47@gatech.edu}.} \and Sam Spiro\footnote{Department of Mathematics, Rutgers University {\tt sas703@scarletmail.rutgers.edu}. This material is based upon work supported by the National Science Foundation Mathematical Sciences Postdoctoral Research Fellowship under Grant No. DMS-2202730.}}

\date{\today}

\theoremstyle{plain}
\newtheorem{thm}{Theorem}[section]
\newtheorem{theorem}[thm]{Theorem}
\newtheorem{conjecture}[thm]{Conjecture}
\newtheorem{lemma}[thm]{Lemma}
\newtheorem{corollary}[thm]{Corollary}
\newtheorem{proposition}[thm]{Proposition}

\theoremstyle{definition}

\newtheorem{question}[thm]{Question}

\newtheorem{remark}[thm]{Remark}

\newtheorem{definition}[thm]{Definition}
\newtheorem{claim}[thm]{Claim}
\newtheorem{fact}[thm]{Fact}

\newtheorem{example}[thm]{Example}

\newtheorem{defn-thm}[thm]{Definition-Theorem}

\newcommand{\btheorem}{\begin{theorem}}
\newcommand{\etheorem}{\end{theorem}}
\newcommand{\bconjecture}{\begin{conjecture}}
\newcommand{\econjecture}{\end{conjecture}}
\newcommand{\bproposition}{\begin{proposition}}
\newcommand{\eproposition}{\end{proposition}}
\newcommand{\bdefinition}{\begin{definition}}
\newcommand{\edefinition}{\end{definition}}
\newcommand{\bcorollary}{\begin{corollary}}
\newcommand{\ecorollary}{\end{corollary}}
\newcommand{\bproof}{\begin{proof}}
\newcommand{\eproof}{\end{proof}}
\newcommand{\bclaim}{\begin{claim}}
\newcommand{\eclaim}{\end{claim}}
\newcommand{\bquestion}{\begin{question}}
\newcommand{\equestion}{\end{question}}
\newcommand{\bfact}{\begin{fact}}
\newcommand{\efact}{\end{fact}}
\newcommand{\bremark}{\begin{remark}}
\newcommand{\eremark}{\end{remark}}
\newcommand{\eexample}{\end{example}}
\newcommand{\bexample}{\begin{example}}

\newcommand{\elemma}{\end{lemma}}
\newcommand{\blemma}{\begin{lemma}}

\newcommand{\beq}{\begin{equation}}
\newcommand{\eeq}{\end{equation}}

\newcommand{\mathify}[1]{\ifmmode{#1}\else\mbox{$#1$}\fi}
\newcommand{\bigO}O

\newcommand{\Y}{\mathcal{Y}}
\newcommand{\C}{\mathcal{C}}

\renewcommand{\l}{\left}
\renewcommand{\r}{\right}
\newcommand{\ex}{\mathrm{ex}}

\newcommand{\p}{{\rm plog}_{\ell,r}(n)}

\newcommand{\Om}{\Omega}

\newcommand{\Ext}{{\rm Ext}}

\newcommand{\remove}[1]{{}}

\begin{document}

\setlength{\baselineskip}{12pt}

\maketitle

\begin{abstract}
An $r$-uniform linear cycle of length $\ell$, denoted by $C^r_{\ell}$, is an $r$-graph with $\ell$ edges $e_1,e_2,\dots,e_{\ell}$ where $e_i=\{v_{(r-1)(i-1)},v_{(r-1)(i-1)+1},\dots,v_{(r-1)i}\}$ (here $v_0=v_{(r-1)\ell}$). For $0<\delta<1$ and $n$ sufficiently large, we show that every $n$-vertex $r$-graph $G$ with $n^{r-\delta}$ edges contains at least $n^{(r-1)(2\ell+1)-\delta(2\ell+1+\frac{4\ell-1}{(r-1)(2\ell+1)-3})-o(1)}$ copies of $C^r_{2\ell+1}$. Further, conditioning on the existence of dense high-girth hypergraphs, we show that there exists $n$-vertex $r$-graphs with $n^{r-\delta}$ edges and at most $n^{(r-1)(2\ell+1)-\delta(2\ell+1+\frac{1}{(r-1)\ell-1})+o(1)}$ copies of $C^r_{2\ell+1}$.
\end{abstract}

\section{Introduction}
Given an $r$-uniform hypergraph ($r$-graph for short) $F$, the Tur\'an number $\ex(F,n)$ is the minimum interger $M$ such that every $n$-vertex $r$-graph with more than $M$ edges contains a subgraph isomorphic to $F$ (referred to as a copy of $F$). Determining or estimating $\ex(F,n)$ is one of the central topics in extremal combinatorics; this kind of problem is known as the {\em Tur\'an problem}. A natural extension of the Tur\'an problem ask: if the number of edges in an $n$-vertex $r$-graph is more than $\ex(n,F)$, then how many copies of $F$ must it necessarily contain? This kind of problem is called the {\em supersaturation problem}.  

For a given $r$-graph $F$, we use $v(F)$ and $e(F)$ to denote its numbers of vertices and edges. Erd\H{o}s and Simonovits~\cite{erdHos1983supersaturated} proved a general conclusion, that is, if an $r$-uniform hypergraph $G$ has $\ex(n,F)+cn^r$ hyperedges, then $G$ contains at least $\varepsilon n^{v(F)}$ copies of $F$ where $\varepsilon$ is depending on $c$, see also~\cite{liu2020exact,reiher2016clique,nikiforov2011number,razborov2008minimal} and references therein. When considering the case of "weakly supersaturated graphs", specifically when $e(G) = \ex(n,F) + o(n^r)$, the scenario becomes more intricate, see~\cite{ma2025supersaturation,mubayi2010counting,pikhurko2017supersaturation} and references therein. In the case of graphs $(r=2)$, Erd\H{o}s and Simonovits made the following conjecture for the supersaturation problem for bipartite graphs:

\begin{conjecture}[\cite{simonovits1984extremal}]\label{conjecture:supersaturation}
Let $F$ be a bipartite graph such that $\ex(n,F)=O(n^{2-\alpha})$ for some $0<\alpha<1$. There exist positive constants $c$ and $c'>0$ such that if $G$ is an $n$-vertex graph such that $e(G)\ge cn^{2-\alpha}$, then $G$ contains at least
\[
c'\frac{e(G)^{e(F)}}{n^{2e(F)-v(F)}}
\]
copies of $F$.
\end{conjecture}

Note that the number $c'\frac{e(G)^{e(F)}}{n^{2e(F)-v(F)}}$ in Conjecture $1.1$ is achieved by the random graph with $e(G)$ edges: it contains in expectation $\binom{n}{v(F)}(\frac{e(G)}{\binom{n}{2}})^{e(F)}\approx c'\frac{e(G)^{e(F)}}{n^{2e(F)-v(F)}}$ copies of $F$. In \cite{faudreecycle}, Simonovits proved Conjecture \ref{conjecture:supersaturation} for even cycles $C_{2\ell}$. In \cite{erdHos1982compactness}, Erd\H{o}s and Simonovits proved it for the path graph $P_4$. Furthermore, in \cite{erdHos1983supersaturated}, they extended their results to many other cases, including the cube graph and the complete bipartite graph $K_{p,q}$. For additional references, the reader is referred to Lov\'{a}sz and Simonovits~\cite{lovasz1983number}, Razborov~\cite{razborov2007flag}, Lov\'{a}sz~\cite{lovasz2012large}, and Reiher~\cite{reiher2016clique}.

It turns out that the supersaturation problem is closely related to the {\em Sidorenko problem}. A \textit{homomorphism} from an $r$-graph $F$ to an $r$-graph $H$ is a map $\phi:V(F)\to V(H)$ such that $\phi(e)$ is an edge of $H$ whenever $e$ is an edge of $F$. Let $\hom(F,H)$ denote the number of homomorphisms from $F$ to $H$ and define the \textit{homomorphism density}
\[t_F(H)=\frac{\hom(F,H)}{v(H)^{v(F)}},\]
which is equivalently the probability that a random map $\phi:V(F)\to V(H)$ is a homomorphism.  From this, we can derive an important definition regarding homomorphism density. 
\begin{definition}
    We say that an $r$-graph $F$ is \textit{Sidorenko} if for every $r$-graph $H$ we have
    \[t_F(H)\ge t_{K_r^r}(H)^{e(F)},\]
    where $K_r^r$ is the $r$-graph consisting of a single edge.
\end{definition}
Observe that if $F$ is not $r$-partite, then $F$ is not Sidorenko (due to $H=K_r^r$, for example). For bipartite graphs, Sidorenko made the following conjecture.

\begin{conjecture}[\cite{sidorenko1993correlation,sidorenko1991inequalities}]
Every bipartite graph is Sidorenko.
\end{conjecture}

A considerable body of research has been devoted to the investigation of Sidorenko's conjecture~\cite{conlon2010approximate,conlon2018some,conlon2017finite,conlon2018sidorenko,coregliano2021biregularity,fox2017local,hatami2010graph,kim2016two,li2011logarithimic,lovasz2011subgraph,szegedy2014information}. Despite these efforts, the conjecture remains widely open.

It is well-established that Sidorenko's conjecture does not generalize to hypergraphs. Indeed, Sidorenko~\cite{sidorenko1993correlation} showed that $r$-uniform linear triangles are not Sidorenko for $r \geq 3$. This result has historically limited interest in identifying which $r$-graphs satisfy the Sidorenko property. However, this trend has recently shifted with the seminal work of Conlon, Lee, and Sidorenko~\cite{conlon2024extremal}, who unveiled a profound connection between non-Sidorenko hypergraphs and Tur\'{a}n numbers.

\begin{thm}[\cite{conlon2024extremal}]\label{thm:CLS}
    If $F$ is not Sidorenko, then there exists $c=c(F)$ such that
    \[\ex(n,F)=\Om\left(n^{r-\frac{v(F)-r}{e(F)-1}+c}\right).\]
\end{thm}
We note that this improves upon the lower bound $\ex(n,F)=\Om(n^{r-\frac{v(F)-r}{e(F)-1}})$ which holds for all $F$ by a standard random deletion argument.

Nie and Spiro~\cite{nie2023sidorenko} introduced the following notion to measure how far an $r$-graph is from being Sidorenko.
\begin{definition}
Given an $r$-partite $r$-graph $F$, let
\[
s(F):=\sup\l\{s: \exists H,~t_F(H)=t_{K_r^r}(H)^{s+e(F)}>0\r\}.
\]
We call $s(F)$ the {\em Sidorenko gap} of $F$.
\end{definition}

Clearly, $F$ is Sidorenko if and only if $s(F)=0$. For non-Sidorenko hypergraphs, Nie and Spiro~\cite{nie2023sidorenko} determined $s(F)$ for $F$ being the $r$-expansion of $k$-uniform clique on $k+1$ vertices (see~\cite{nie2023sidorenko} for precise definition). 

One can show that if $F$ is not Sidorenko, then the analogue of~\Cref{conjecture:supersaturation} for $F$ is not true.
More precisely, we can prove the following.
\begin{proposition}\label{Copies of F}
For every $r\geq 3$, let $F$ be an $r$-graph. For any integer $n$ and reals $\varepsilon, \delta>0$, there exists $N>n$ such that there exists an $N$-vertex $r$-graphs $H_N$ with at least $N^{r-\delta}$ edges such that $H_N$ contains at most
\[
N^{v(F)-\delta(s(F)+e(F)-\varepsilon)}
\]
copies of $F.$
\end{proposition}

Inspired by \Cref{Copies of F}, we make the following conjecture.
\begin{conjecture}\label{conjecture:supersat}
For every $r\geq 3$, let $F$ be an $r$-partite $r$-graph such that $\ex(F,n)=O(n^{r-\alpha})$ for some $0<\alpha<r-1$. Then there exists positive constant $c$ such that if $G$ is an $n$-vertex $r$-graph such that $e(G)\ge cn^{\alpha}$, then $G$ contains at least
\[
n^{v(F)-o(1)}\l(\frac{e(G)}{n^r}\r)^{e(F)+s(F)}
\]
copies of $F$.
\end{conjecture}

An $r$-uniform linear cycle of length $\ell$, denoted by $C^r_{\ell}$, is an $r$-graph with $\ell$ edges $e_1,e_2,\dots,e_{\ell}$ where $e_i=\{v_{(r-1)(i-1)},v_{(r-1)(i-1)+1},\dots,v_{(r-1)i}\}$ (here $v_0=v_{(r-1)\ell}$). \Cref{conjecture:supersat} have been implicitly confirmed for $C^r_{2\ell},~\ell\ge 2$~\cite{nie2024turan,mubayi2023random}, and also for the non-Sidorenko hypergraph the $r$-expansion of $k$-uniform clique on $k+1$ vertice~\cite{balogh2019number}~\cite{nie2023random}.

In this paper, we study the supersaturation problem and the Sidorenko gap problem of odd linear cycles. The problems for $C^r_3$ have been solved in~\cite{nie2021triangle,nie2023random,nie2023sidorenko}; it is proved that $s(C^r_3)=\frac{1}{r-2}$ and \Cref{conjecture:supersat} has been confirmed for $C^r_3$. The problems for $C^r_{2\ell+1}$ where $\ell\ge 2$ are still wide open. In~\cite{nie2023sidorenko}, Nie and Spiro showed that $s(C^r_{2\ell+1})\le \frac{2\ell-1}{r-2}$. Further, they mentioned (without a proof) a conditional lower bound $s(C^r_{2\ell+1})\ge \frac{1}{(r-1)\ell-1}$. Here we prove this lower bound formally and as a consequence obtain a conditional upper bound for the supersaturation problem. In this paper $o(1)$ means a function that goes to $0$ when $n$ tends to infinity. A Berge-Cycle of length $k$ in a hypergraph is a set of $k$ distinct vertices $\{v_1,v_2,\dots,v_k\}$ and $k$ distinct edges $\{e_1,e_2,\dots,e_k\}$ such that $\{v_i,v_{i+1}\} \subset e_i$ with indices
 taken modulo $k$. We call the length of the shortest Berge-Cycle contained in a hypergraph the girth of this hypergraph.

\begin{thm}\label{theorem:main_lowerbound}
For every $r\geq 3$ and $\ell\geq 2$, if for all positive integers $n$ there exist $n$-vertex $r$-graphs of girth $2\ell+2$ with $n^{1+1/\ell-o(1)}$ edges, then
\[
s(C^r_{2\ell+1})\ge \frac{1}{(r-1)\ell-1},
\]
and thus, for every constant $0<\delta<1$, there exists $n$-vertex $r$-graphs with at least $n^{r-\delta}$ edges containing at most
\[
n^{(r-1)(2\ell+1)-\delta(2\ell+1+\frac{1}{(r-1)\ell-1})+o(1)}
\]
copies of $C^r_{2\ell+1}$.
\end{thm}
The condition of the existence of dense high-girth hypergraphs here if true would be best possible~\cite{collier2018linear}. It can be viewed as a hypergraph analogue of the Erd\H{o}s-Simonovits conjecture for the existence of dense high-girth graphs, which is widely believed to be true among the research community. In the case when $\ell=1$, such condition is true due to the seminal Ruzsa-Szemer\'edi construction~\cite{ruzsatriple} and its generalization~\cite{gowers2021generalizations}. See also \cite{lazebnik2003hypergraphs, shangguan2020sparse, haymaker2024hypergraphs} for more on this topic.

On the other hand, concerning lower bounds for the supersaturation problem, as a consequence of~\cite[Theorem 2.5]{nie2023random}, we know that an $r$-graph with $e(G)\ge n^{r-1+o(1)}$ contains at least 
\[
n^{(r-1)(2\ell+1)-o(1)}\l(\frac{e(G)}{n^r}\r)^{2\ell+1+\frac{2\ell-1}{r-2}}
\]
copies of $C^r_{2\ell+1}$. We prove the following better lower bounds for $\ell\ge 2$.

\begin{thm}\label{theorem:main_upperbound}
For every $r\geq3$, $\ell\geq 2$, there exists a constant $C$ such that if $G$ is an $n$-vertex $r$-graph such that $e(G)\ge n^{r-1}(\log n)^{C}$, then, as $n\rightarrow\infty$, $G$ contains at least
\[
n^{(r-1)(2\ell+1)-o(1)}\l(\frac{e(G)}{n^{r}}\r)^{\l(2\ell+1+\frac{4\ell-1}{(r-1)(2\ell+1)-3}\r)}
\]
copies of $C^r_{2\ell+1}$. As a consequence,
we have 
\[
s(C^r_{2\ell+1})\le \frac{4\ell-1}{(r-1)(2\ell+1)-3}.
\]
\end{thm}

Theorem~\ref{theorem:main_upperbound} is our main result. In fact, it is obtained by enhancing a weaker result using the method of codegree dichotomy. We present this weaker result here.

\begin{proposition}\label{theorem:weaker main_upperbound}
For every $r\geq3$, $\ell\geq 2$, there exists a constant $C$ such that if $G$ is an $n$-vertex $r$-graph such that $e(G)\ge n^{r-1}(\log n)^C$, then, as $n\rightarrow\infty$, $G$ contains at least
\[
n^{(r-1)(2\ell+1)-o(1)}\l(\frac{e(G)}{n^r}\r)^{(2\ell+1+\frac{2}{r-1})}
\]
copies of $C^r_{2\ell+1}$. As a consequence,
we have 
\[
s(C^r_{2\ell+1})\le \frac{2}{r-1}.
\]
\end{proposition}

It is worth noting that Theorem~\ref{theorem:main_upperbound} and Proposition~\ref{theorem:weaker main_upperbound} give the same bound when $r=3$.

\subsection{Notation}

A $k$-shadow of $H$ is a set $X$ of $k$ vertices such that $X\subset e$ for some
edge $e$ of $H$. The codegree of $X$ is the number of edges of $H$ containing $X$, which is written as $d_H(X)$; when $X=\{u, v\}$, we simplify the notation to $d_{H}(u, v)$. We denote by $\Delta_j(H)$ the maximum $d_H(\delta)$ among all $j$-shadows $\delta$. Whenever $H$ is clear from the context we will drop it from the notation. For an $r$-partite $r$-graph $G$ with vertex partition $(V_1, V_2, \cdots , V_r)$, for $1\leq i<j\leq r$, we write $\partial_{V_i, V_j}(H)$ for the set of $2$-shadows $(v_i, v_j)$ with $v_i\in V_i$ and $v_j\in V_j$. Usually we simplify the notation to $\partial_{i,j}(H)=\partial_{V_i, V_j}(H)$. Throughout this paper, We will use a notation $\p$ to denote any function which is of order $\Theta((\log n)^{f(\ell,r)})$.

For an $r$-graph $H$ and a positive integer $t$, the $t$-$blow$-$up$ $H[t]$ is the $r$-graph on vertex set $V(H)\times [t]$ with hyperedges $\{\{(v_1,t_1), \dots , (v_r, t_r)\}: \{v_1, \dots, v_r\}\in H, 1\leq t_1, \dots, t_r \leq t \}$.

\subsection{Related work}

Fox, Sah, Sawhney, Stoner, and Zhao~\cite{fox2020triforce} studied the supersaturation problem of $C^3_3$, but in a slightly different flavor: they showed that for any $t>0$, a sufficiently large 3-graph $G$ with edge density $t$ must contain at least 
$|V(G)|^4t^{4-o(1)}$ copies of $C^3_3$.

Jiang and Yepremyan~\cite{jiang2020supersaturation} studied the supersaturation of even linear cycles in linear hypergraphs. Recently, Mubayi and Solymosi studied the supersaturation of linear 5-cycles in linear hypergraphs~\cite{mubayi2025many}.
\section{Proof}

The proof idea for Proposition~\ref{Copies of F} is to use the blow-up technique to link the homomorphic copies of $F$ in $H$ with the copies of $F$ in $H[t]$, and then employ the Tensor product of hypergraphs to ensure that the hypergraph can be sufficiently large.

\subsection{Proof of Proposition~\ref{Copies of F}}

\begin{proof}[Proof of Proposition~\ref{Copies of F}]
    
By the definition of $s(F)$, for any $\varepsilon>0$, there exists some $\delta'>0$ and an $m$-vertex $r$-graph $H$ with $m^{r-\delta'}$ edges such that 
\[
\hom(F,H)\leq m^{v(F)-\delta'(e(F)+s(F)-\frac{\varepsilon}{2})}.
\]

Moreover, we define the tensor product $H\otimes H=H^2$ to be $r$-graph on $V(H)\times V(H)$, where $((x_1,y_1),\dots , (x_r,y_r))\in E(H^2)$ if and only if $(x_1,\dots, x_r),\ (y_1,\dots ,y_r)\in E(H)$.
Clearly we have $v(H^2)=v(H)^2=m^2$ and $e(H^2)= e(H)^2=m^{2(r-\delta')}$. 
Furthermore, each copy of $F$ in $H^2$ induces two (possibly the same) copies of $F$ in $H$. 
Hence, we have $\hom(F,H^2)\leq \hom(F,H)^2 = m^{2(v(F)-\delta'(s(F)+e(F)-\frac{\varepsilon}{2}))}$, which means $H^2$ still satisfies the upper bound of the number of $F$-copies. 
Therefore, by repeatedly applying the tensor product operation, we may assume that $m$ is sufficiently large.

Now we construct the desired $r$-graph promised by the proposition based on $H$.
If $\delta\leq\delta'$, then 
let $H[t]$ denote the $t$-$blow$-$up$ of $H$, where the vertex set is given by \[\{v_{ij}, i\in[m], j\in [t]\}\] and let $\phi$ be an arbitrary homomorphism from $F$ to $H$. 
We define an injection $\phi'$ from $F$ to $H[t]$, if there exists $v_a, v_b\in V(F)$, such that $\phi(v_a)=\phi(v_b)=v_i\in V(H)$,  then we assign $\phi'(v_1)$ and $ \phi'(v_2)$ to different vertices within the same $t$-set. 
i.e. $\phi'(v_1)=v_{ij}\in V(H[t])$, $ \phi'(v_2)=v_{ik}\in V(H[t])$, where $1\leq j\neq k\leq t$.  
Clearly, such a $\phi'$ corresponds to a copy of $F$ in $H[t]$. Conversely, each labeled copy of $F$ in $H[t]$ also corresponds to such a homomorphism. 
When $t$ is sufficiently large, there are at most $\hom(F,H)\cdot t^{v(F)}$ such injections. 
Moreover, let $n=mt$ and $\delta''=\frac{\log m}{\log n}\delta'$. 
Note that $m^{\delta'}=n^{\delta''}$. Clearly we have that $H[t]$ has $n$ vertices and $m^{r-\delta'}\cdot t^r=n^{r-\delta''}$ edges. 
We find that $H[t]$  contains at most $m^{v(F)-\delta'(e(F)+s(F)-\frac{\varepsilon}{2})}\cdot t^{v(F)}\leq n^{v(F)-\delta''(e(F)+s(F)-\varepsilon)}$ labeled copies of $F$. 
Setting $t=m^{\frac{\delta'}{\delta}-1}$ yields $\delta''=\delta$. 
Thus, we obtain the desired result.

If $\delta>\delta'$, then we consider $H[p]$, a random induced subgraph of $H$, where $p$ will be specified later.
That is, $H[p]$ is obtained from $H$ by keeping each vertex in $H$ is independently and randomly with probability $p$. 
Set $n=mp$ and let $Y$ be the number of $F$-copies in $H[p]$.
Let $\delta''=\frac{\delta'\log m+2r\log\frac{3}{2}}{\log n}$ and note that $m^{-\delta'}=(\frac{3}{2})^{2r}n^{-\delta''}$. 
Clearly we have $\mathbb E(v(H[p]))=mp=n$, $\mathbb E(e(H[p]))=m^{r-\delta'}p^r=(\frac{3}{2})^{2r}n^{r-\delta''}$ and 
\[
\begin{aligned}
\mathbb E(Y)&\leq \hom(F,H)\cdot p^{v(F)}\\
&\leq m^{v(F)-\delta'(e(F)+s(F)-\frac{\varepsilon}{2})}\cdot p^{v(F)}\\
&=(3/2)^{2r(e(F)+s(F)-\frac{\varepsilon}{2})}n^{v(F)-\delta''(e(F)+s(F)-\frac{\varepsilon}{2})}.    
\end{aligned}
\]

Since the number of vertices $N$ in $H[p]$ follows a binomial distribution $N \sim \text{Bin}(m, p)$, by Chernoff's Bound, we have 
\[
\Pr\left(|N - n| \geq \frac{n}{2}\right) \leq 2\exp\left(-\frac{n}{12}\right).
\]
Similarly as $e(H[p]) \sim \text{Bin}(m^{r-\delta'}, p^r)$, by Chernoff's Bound, we have 
\[
\Pr\left(e(H[p]) < \frac{1}{2}\mathbb E[e(H[p])]\right) \leq \exp\left(-\frac{\mathbb E[e(H[p])]}{8}\right).
\]
Moreover, for $Y$, by Markov's inequality, we have
\[
\Pr(Y \geq 2\mathbb E[Y]) \leq \frac{1}{2}.
\]
Setting $p=m^{\frac{\delta'}{\delta}-1}\cdot (\frac{3}{2})^{\frac{2r}{\delta}}$ yields $\delta''=\delta$. 
Let $A,B,C$ be the events that $N\in [\frac{n}{2},\frac{3n}{2}]$, $e(H[p]) \geq \frac{1}{2}\mathbb E[e(H[p])]$, and $Y < 2\mathbb E[Y]$. 
We find that as $m$ is sufficiently large, since $n=mp=m^{\frac{\delta'}{\delta}}\cdot (\frac{3}{2})^{\frac{2r}{\delta}}$, the following holds:
$\mathbb E(v(H[p]))\geq 24\ln 8$ and $\mathbb E(e(H[p]))\geq 8\ln 8$.
Thus, we infer that each of $A$ and $B$ occurs with probability at least $7/8$.
Hence, with probability at least $1/4$, events $A, B$ and $C$ hold for $H[p]$.
This yields the existence of the desired subgraph $H[p]$ with $e(H[p])\geq \frac{1}{2}(\frac{3}{2})^{2r}n^{r-\delta''}\geq (\frac{3}{2}n)^{r-\delta}$ which contains at most 
\[
2\left(\frac{3}{2}\right)^{2r(e(F)+s(F)-\frac{\varepsilon}{2})}n^{v(F)-\delta''(e(F)+s(F)-\frac{\varepsilon}{2})}\leq \left(\frac{n}{2}\right)^{v(F)-\delta(e(F)+s(F)-\varepsilon)}
\]
copies of $F$. 
%
%
%
\end{proof}

\subsection{Proof of Theorem~\ref{theorem:main_lowerbound}}

Before presenting the proof, we first provide some necessary definitions.
A path of a linear hypergraph is a sequence $v_1e_1v_2e_2\cdots v_{k-1}e_{k-1}v_k$, where for any $1\leq i\leq k-1$, $\{v_i,v_{i+1}\}\subset e_i$, and each $v_i$ appears only once in this sequence. The length of a path is the number of edges of the path.
We say that two vertices $v_1$ and $v_2$ in a linear hypergraph is connected if there exists a path containing $v_1$ and $v_2$, and we say that a linear hypergraph $H$ is connected if any two vertices in $H$ are connected. We say that $v_1$ is at distance $i$ from $v_2$ if the length of the shortest path containing $v_1$ and $v_2$ is $i$.

Firstly, we present two simple propositions.

\begin{proposition}\label{proposition:homo}
 Let $H$ be an $n$-vertex $r$-graph of girth $2\ell+2$ with $n^{1+1/\ell-o(1)}$ edges. If $\phi$ is a homomorphism from $C^r_{2\ell+1}$ to $H$, then $\phi(V(C^r_{2\ell+1}))$ induces a linear tree in $H$ with at most $\ell$ edges.
\end{proposition}

\begin{proof}[Proof of Proposition~\ref{proposition:homo}]
    Note that $H$ is a linear hypergraph since we regard two edges with an intersection of at least $2$ vertices as forming a Berge-Cycle of length $2$. Moreover, $\phi(V(C^r_{2\ell+1}))$ induces a linear tree in $H$ since the girth of $H$ is $2\ell+2$, and if $\phi(V(C^r_{2\ell+1}))$ forms a cycle, its length must be at most $2\ell+1$.

     Let us prove that the number of edges of $\phi(V(C^r_{2\ell+1}))$ is at most $\ell$. 
    Let the linear tree induced by $\phi(V(C^r_{2\ell+1}))$ be denoted as $T$. Suppose $e(T)\geq \ell+1$. Since $\phi$ is a homomorphism and $e(C^r_{2\ell+1})=2\ell+1$, by the pigeonhole principle, there must be at least one edge in $E(T)$ whose preimage consists of exactly one edge, say $e_1$. 
    We remove $e_1$ from $C^r_{2\ell+1}$ and denote the resulting hypergraph as $C^-$. 
    Let the two adjacent edges of $e_1$ be denoted as $e_2$ and $e_3$, with $e_1\cap e_2=v_1$ and $e_1\cap e_3=v_2$. Since $\{v_1,v_2\} \subset e_1$, $\phi(v_1)$ and $\phi(v_2)$ are distinct vertices of $\phi(V(e_1))$ in $T$. Therefore, if we remove $\phi(V(e_1)\setminus\{v_1,v_2\})$ from $T$, $v_1$ and $v_2$ will become disconnected, as $T$ is a linear tree. 
    Furthermore, observe that $C^-$ remains connected. Clearly, the homomorphic image of a connected hypergraph is still connected, which leads to a contradiction. Hence, we conclude that $e(T)\leq \ell$.
    \end{proof}

\begin{proposition}\label{proposition:countingtree}
There exists $H'$ with $n^{1-o(1)}$ vertices, $n^{1+1/\ell-o(1)}$ edges, and the number of linear trees in $H'$ with at most $\ell$ edges is at most $n^{2+o(1)}.$
\end{proposition}

 \begin{proof}[Proof of Proposition~\ref{proposition:countingtree}]
    We iteratively remove vertices from $H$ whose degree is less than $\frac{1}{2}n^{1/\ell-o(1)}$, until the degree of all remaining vertices is greater than $\frac{1}{2}n^{1/\ell-o(1)}$.
    During the process, at most $n\cdot \frac{1}{2}n^{1/\ell-o(1)}<\frac{1}{2}n^{1+1/\ell-o(1)}$ edges are removed. 
    Thus, we can obtain a subhypergraph $H'\subset H$ with at least $\frac{1}{2}n^{1+1/\ell-o(1)}$ edges, and the girth of $H'$ remains at least $2\ell+2$. Additionally, the minimum degree of $H'$ is $\delta(H')\geq \frac{1}{2}n^{1/\ell-o(1)}$.

      Suppose $x\in V(H')$ and $d(x)=\Delta(H')$. Let $U_i$ be the set of vertices that are at distance $i$ from $x$.
    Let $m_i$ be the number of vertices in $U_i$. Clearly, $m_1=\Delta(H')$.

      Let $d_{[U,V]}(v)$ denote the number of edges in $H'$ that are restricted to $U$ and $V$ and contain the vertex $v$. Note that for any $v_{i+1}\in U_{i+1}, 1\leq i\leq \ell-1$, $d_{[U_i,U_{i+1}]}(v_{i+1})=1$. Otherwise, there would exist vertices $v_a, v_b\in U_i$ and edges $e_1, e_2\in E(H')$, such that $\{v_a,v_{i+1}\}\subset e_1$ and $\{v_b, v_{i+1}\}\subset e_2$. This would imply the existence of a cycle containing $e_1, e_2$ with length less than $2\ell+2$, which is a contradiction. 

      Moreover, $e(H'[U_i])=0$, $1\leq i\leq \ell$. Otherwise, there exists an edge $e\in E(H'[U_i])$, then there would exist a cycle containing $e$ with length less than $2\ell+2$, which is a contradiction.

      Because for any $v_i\in U_i, d(v_i)\geq \delta(H')\geq \frac{1}{2}n^{1/\ell-o(1)}$ and $d(v_i)=d_{[U_i,U_{i-1}]}(v_i)+d_{H'[U_i]}(v_i)+d_{[U_i,U_{i+1}]}(v_i)$, we have $d_{[U_i,U_{i+1}]}(v_i)\geq \frac{1}{2}n^{1/\ell-o(1)}-1$. Thus, we have 
    \begin{equation}\label{eq:V_i}   
    m_{i+1}\geq \frac{1}{3}n^{1/\ell-o(1)}m_i,\  \mbox{for any}\  1\leq i\leq \ell-1.
     \end{equation}

       Thus, we have $n\geq m_\ell\geq \Delta(H')\cdot (\frac{1}{3}n^{1/\ell-o(1)})^{\ell-1}=(\frac{1}{3})^{\ell-1}n^{1-1/\ell-o(1)}\cdot\Delta(H')$, which means that \[\Delta(H')\leq n^{1/\ell+o(1)}.\]

     Moreover, we have $m_{\ell}\geq (\frac{1}{3})^{\ell-1}n^{1-1/\ell-o(1)}\cdot\Delta(H')\geq (\frac{1}{3})^{\ell-1}n^{1-1/\ell-o(1)}\cdot\delta(H')\geq (\frac{1}{3})^{\ell-1}\cdot \frac{1}{2}n^{1-o(1)}$, which implies that $|V(H')|=n^{1-o(1)}$.

       In conclusion, the number of linear trees with at most $\ell$ edges is at most 
     \[
     n\cdot \Delta(H')^\ell\cdot f(\ell,r)\leq n^{2+o(1)}.
     \]
     where $f(\ell,r)$ is a function bounded by $\ell$ and $r$, representing the number of choices for selecting the vertex other than the first one. 
    \end{proof}

\begin{proof}[Proof of Theorem~\ref{theorem:main_lowerbound}]
     Let $H$ be an $n$-vertex $r$-graph of girth $2\ell+2$ with $n^{1+1/\ell-o(1)}$ edges as assumed in the theorem.

      By the definition of the homomorphism density,  proposition~\ref{proposition:homo} and proposition~\ref{proposition:countingtree}, we have 
    \[
    t_{C_{2\ell+1}^r}(H')=\frac{\hom(C_{2\ell+1}^r,H')}{v(H')^{v(C_{2\ell+1}^r)}}=\frac{g(\ell,r)n^{2+o(1)}}{n^{(2\ell+1)(r-1)-o(1)}}.
    \]
    where $g(\ell,r)$ is a function bounded by $\ell$ and $r$. And 
    \[
    t_{K^r_r}(H')=\frac{r!n^{1+1/\ell-o(1)}}{n^{r-o(1)}}.
    \]
       By definition, $s(C^r_{2\ell+1})$ is at least $\log(t_{C_{2\ell+1}^r}(H'))/\log(t_{K^r_r}(H'))-e(C_{2\ell+1}^r)$. Thus, we have 
     \[
     s(C^r_{2\ell+1}) \geq \frac{\log \phi(\ell,r)+(2+o(1))\log n- ((2\ell+1)(r-1)-o(1))\log n}{\log r!+(1+1/\ell-o(1))\log n-(r-o(1))\log n}-(2\ell+1),
     \]
     Taking $n\rightarrow \infty$, we have
     \[
     s(C^r_{2\ell+1}) \geq \frac{2-(2\ell+1)(r-1)}{1+1/\ell-r}-(2\ell+1)=\frac{1}{(r-1)\ell-1}.
     \]
\end{proof}

\subsection{Proof of Proposition~\ref{theorem:weaker main_upperbound}}

Our proof for Proposition~\ref{theorem:weaker main_upperbound} primarily draws upon the methodology employed by Mubayi and Yepremyan~\cite{mubayi2023random} in their treatment of balanced supersaturation for even cycles in $r$-graphs. We have extended this approach to address the issue of supersaturation for odd cycles in $r$-graphs.

In our proof, we need a supersaturation result for even cycles in simple graphs by Simonovits~\cite{erdos1984cube} (unpublished, see \cite[Theorem 1.5]{morris2016number} for a proof).

\begin{thm}[\cite{erdos1984cube,morris2016number}]\label{thm:morrissaxton}
For every $\ell\geq 2$, there exist $\delta >0$ and $k_0\in \mathbb{N}$ such that   the following holds for every $k\geq k_0$ and every $n\in \mathbb{N}$. Given a graph $G$ with $n$ vertices and $kn^{1+1/\ell}$ edges there exists a collection $\mathcal{F}$  of copies of $C_{2\ell}$ in $G$ such that $|\mathcal{F}|\geq \delta k^{2\ell}n^2$.

\end{thm}

For convenience, we restate \Cref{theorem:weaker main_upperbound}.

\begingroup
\def\thethm{\ref{theorem:weaker main_upperbound}}
\begin{proposition}
For every $r\geq3$, $\ell\geq 2$, there exists a constant $C$ such that if $G$ is an $n$-vertex $r$-graph such that $e(G)\ge n^{r-1}(\log n)^C$, then, as $n\rightarrow\infty$, $G$ contains at least
\[
n^{(r-1)(2\ell+1)-o(1)}\l(\frac{e(G)}{n^r}\r)^{(2\ell+1+\frac{2}{r-1})}
\]
copies of $C^r_{2\ell+1}$. As a consequence,
we have 
\[
s(C^r_{2\ell+1})\le \frac{2}{r-1}.
\]
\end{proposition}
\endgroup

\begin{proof}[Proof of \Cref{theorem:weaker main_upperbound}]
    Let  $Q, \delta_0, k_0$ be the quantities derived from Theorem~\ref{thm:morrissaxton} when applied with $\ell$. Let $G$ be a given $r$-graph with $n$ vertices such that $e(G)\geq n^{r-1}(\log n)^C$, where $C$ is a constant. According to a classical result by Erd\H{o}s and Kleitman~\cite{erdos1968coloring}, $G$ contains an $r$-partite subgraph $H$ with $r$-partition $(V_1,V_2,\dots, V_r)$
such that the number of edges in $H$ satisfies $e(H)\geq r!e(G)/r^r$.

Let $\Y$ be the collection of all functions from $\binom{[r]}{2}$ to $[\log \binom{n}{r-2}]$, and for each $i,j\in [r]$ and $y\in\Y$ we simply write $y_{ij}$ for $y(i,j)$. We say that an edge $e\in H$ is of $type$ $y$ if for every set $\{i,j\}\subseteq \binom{[r]}{2}$, we have $2^{y_{ij}-1}\leq d_{H}(e\cap (V_i\cup V_j))<2^{y_{ij}}$. By the pigeonhole principle, there exists some $y\in \Y$ such that at least $\log (n^r)^{-r^2}e(H)$ edges of $H$ are of $type$ $y$. Let $H_0$ be the collection of edges of $type$ $y$, and let  $\Delta_{ij}\coloneqq 2^{y_{ij}}$. 

Set $H^*=H_0$ and let $U_i\coloneqq V_i\cap V(H^*)$. If there exists some pairs $v_i\in U_i$,  $v_j\in U_j$, $1\leq i<j\leq r$ such that $0<d_{H_0}(v_i,v_j)<2^{-r^2}\log (n^r)^{-r^2}d_H(v_i,v_j)$ , then delete every edge from $H^*$ containing $v_i$ and $v_j$. We excute this operation iteratively on $H^*$. After this, we also remove all the vertices which became isolated in the resulting hypergrah. Let $H'$ be the resulting hypergraph, we continue to use the notation $U_i$ for the parts of $H'$. Noting that it satisfies that if  $v_i\in V_i, v_j\in V_j$  such that there exists $e\in E(H')$ with $\{v_iv_j\}\subseteq e$  then  
\begin{equation}\label{bound of codegree}
    \frac{\Delta_{ij}}{2(2\log{n^r})^{r^2}}\leq \frac{d_H(v_i,v_j)}{(2\log n^r)^{r^2}}\leq d_{H'}(v_i,v_j) \leq \Delta_{ij}.
\end{equation}
Furthermore, note that the total number of edges deleted from $H_0$ to obtain $H'$ is at most 
\[\sum_{1\leq i<j\leq r}\sum_{v_i\in V_i, v_j\in V_j}{\frac{d_H(v_i,v_j)}{(2\log{n^r})^{r^2}}} \leq \frac{r^2e(H)}{(2\log{n^r})^{r^2}} \leq \frac{e(H_0)}{2}.\] 
Combining the above, we have:
\begin{equation}\label{lowerbound of edges}
    e(H')\geq \frac{e(H_0)}{2}\geq \frac{e(H)}{2\log (n^r)^{r^2}}\geq  \frac{e(G)}{\p}.
\end{equation}
Let $\partial_{ij}= \partial_{U_i, U_j}(H')$, without loss of generality, we may assume that $|U_1|\geq |U_2|\geq \cdots \geq|U_r|$.

\begin{claim}\label{the third vertex} For arbitrary $u_1\in U_1$, $u_2\in U_2$, there exists $3\leq i\leq r$ such that there are at least $d_{H'}(u_1,u_2)^{1/(r-2)}$ distinct $x\in U_i$, the set $\{u_1,u_2,x\}$ is part of some edges in $H'$.
\end{claim}

\begin{proof}
Assume that for all $3\leq i\leq r$, we can find fewer than $d_{H'}(u_1,u_2)^{1/(r-2)}$ distinct $x\in U_i$ that satisfy the given conditions. Thus, the number of edges containing $u_1$ and $u_2$ is less than \[(d_{H'}(u_1,u_2)^{1/(r-2)})^{r-2}=d_{H'}(u_1,u_2),\]
which is a contradiction.
\end{proof}

\begin{claim} \label{keyclaim} There exists some $j\in\{2,\dots, r\}$ such that \[|\partial_{1j}| \geq \frac{e(G)^{\frac{1}{r-1}} |U_1|^{\frac{r-2}{r-1}}}{\p}.\]
\end{claim}

\begin{proof}
    Assume that $j=2$ does not meet the desired bound, as otherwise the proof is trivial. We aim to demonstrate the existence of some $j\in \{3, \dots, r\}$ that can be made large enough.

    By (\ref{bound of codegree}), (\ref{lowerbound of edges}), and the definition of $\partial_{12}$, we have:
    \begin{equation}
        \frac{e(G)}{\p}\leq e(H')\leq |\partial_{12}|\Delta_{12}.
    \end{equation}

    Consequently,
    \begin{equation}\label{bound of Delta12}
        |\Delta_{12}|\geq \frac{e(G)}{|\partial_{12}|\p}.    
    \end{equation}
    Given $u_1\in U_1, u_2\in U_2$, if $d_{H'}(u_1,u_2)>0$, then combining (\ref{bound of codegree}) and (\ref{bound of Delta12}), we arrive at:
    \[d_{H'}(u_1,u_2)\geq \frac{\Delta_{12}}{2(2\log{n^r})^{r^2}} \geq \frac{e(G)}{2(2\log{n^r})^{r^2}|\partial_{12}|\p}=\frac{e(G)}{|\partial_{12}|\p}.\]
    Note that $u_1$ and $u_2$ are arbitrary vertices in $U_1$ and $U_2$, respectively, and they are part of at least one edge. Each such edge takes the form $\{u_1,u_2,u_3,\dots, u_r\}$, with $u_j\in U_j$ for $j=3,\dots, r$. Claim~\ref{the third vertex} guarantees the existence of $j$ such that for at least $d_{H'}(u_1,u_2)^{1/(r-2)}$ distinct $x\in U_j$, the set $\{u_1, u_2,x\}$ is a $3$-shadow. Therefore, for the chosen $u_1\in U_1$, we identify some $j=j(u_1)\in \{3,\dots, r\}$ and a subset $X_{j,u_1}\subseteq U_j$ such that $u_1, x_j$ are in an edge for each $x_j\in X_{j,u_1}$, and 
    \[|X_{j,u_1}|\geq d_{H'}(u_1,u_2)^{1/(r-2)}\geq \left(\frac{e(G)}{\p|\partial_{12}|}\right)^{1/(r-2)}.\]
    
    Applying the pigeonhole principle, there exists $j\in\{3,\dots, r\}$ such that for at least $|U_1|/r$ vertices $u_1\in U_1$, $j(u_1)=j$. Observe that for each $u_1\in U_1$ and $x\in X_{j,u_1}$, we obtain distinct edges $u_1x\in \partial_{1j}$. By assumption:
    \[|\partial_{12}| \leq \frac{e(G)^{\frac{1}{r-1}} |U_1|^{\frac{r-2}{r-1}}}{\p},\]
    we have

\[|\partial_{1j}|\geq \frac{|U_1|}{r}\left(\frac{e(G)}{\p|\partial_{12}|}\right)^{1/(r-2)} \geq   \frac{|U_1|}{r}\left(\frac{e(G)}{\p|U_1|^{\frac{r-2}{r-1}} e(G)^{\frac{1}{r-1}}}\right)^{1/(r-2)} = \frac{e(G)^{\frac{1}{r-1}} |U_1|^{\frac{r-2}{r-1}}}{\p},\]
as desired.
\end{proof}

  By Claim~\ref{keyclaim}, we may assume  without loss of generality that

\begin{equation}\label{eq:partialedges}|\partial_{12}| \geq  \frac{e(G)^{\frac{1}{r-1}} |U_1|^{\frac{r-2}{r-1}}}{\p} \geq \frac{e(G)^{\frac{1}{r-1}} m^{\frac{r-2}{r-1}}}{\p},
\end{equation} 
where $m=|U_1|+|U_2|$.

Set $k=\frac{|\partial_{12}|}{m^{1+1/\ell}}$. By (\ref{eq:partialedges}) and $\ell\geq 2, r\geq 3$, we can ensure that 
$k \geq k_0,$  thus by Theorem~\ref{thm:morrissaxton} applied with $k$ and $m$, the  shadow graph $\partial_{12}$ contains a collection $\mathcal{F}$  of copies of $C_{2\ell}$ satisfying $|\mathcal{F}|\geq \delta_0 k^{2\ell}m^2$.


Let $C$ be any $2\ell$-cycle in  $\mathcal{F}$ with consecutive edges $x_1x_2\dots x_{2\ell}$ in the natural cyclic order. i.e., it contains edges $\{x_1,x_2\}, \{x_2,x_3\}, ..., \{x_{2\ell},x_1\}$. By (\ref{bound of codegree}), which states $\frac{\Delta_{12}}{2(2\log{n^r})^{r^2}}\leq d_{H'}(v_1,v_2) \leq \Delta_{12}$, there exists at least $\frac{\Delta_{12}}{2(2\log{n^r})^{r^2}}$ different hyperedges containing both $x_1$ and $x_{2\ell}$. By Claim~\ref{the third vertex}, there exists $3\leq a\leq r$ such that there are at least $(\frac{\Delta_{12}}{2(2\log{n^r})^{r^2}})^{1/(r-2)}$ different $x_{2\ell+1}\in U_a$ for which there exists an edge $e\in E(H')$ containing both $ \{x_{2\ell}, x_{2\ell+1}\} $ and $ \{x_{2\ell+1}, x_1\} $. Therefore, every cycle $C$ can be extended to $(\frac{\Delta_{12}}{2(2\log{n^r})^{r^2}})^{1/(r-2)}$ different sets of the form $\{\{x_1,x_2\}, \{x_2,x_3\}, ..., \{x_{2\ell},x_{2\ell+1}\}, \{x_{2\ell+1}, x_1\}\}$. Let the collection of such sets be denoted by $\C$.

  For every $C'\in\C$, since  $d_{H'}(x_i,x_{i+1})\geq \frac{\Delta_{12}}{2(2\log{n^r})^{r^2}}$ for every $1\leq i \leq 2\ell-1$,   $d_{H'}(x_1,x_{2\ell+1})\geq \frac{\Delta_{1a}}{2(2\log{n^r})^{r^2}}$ and $d_{H'}(x_2,x_{2\ell+1})\geq \frac{\Delta_{2a}}{2(2\log{n^r})^{r^2}}$, the number of ways to extend $C'$ to some linear $(2\ell+1)$-cycles in $H'$ is at least

\begin{equation}\label{extend}
\left(\frac{\Delta_{12}}{2(2\log{n^r})^{r^2}}-2\ell rn^{r-3}\right)^{2\ell-1} \left(\frac{\Delta_{1a}}{2(2\log{n^r})^{r^2}}-2\ell rn^{r-3}\right)\left(\frac{\Delta_{2a}}{2(2\log{n^r})^{r^2}}-2\ell rn^{r-3}\right).
\end{equation}

Indeed, as it suffices to choose all these new vertices to be distinct, each time when we pick a new edge of the cycle we must avoid at most $2\ell(r-1)+1\le 2\ell r$ vertices, which makes at most $2\ell r n^{r-3}$ $(r-2)$-sets unavailable.
Furthermore, by (\ref{bound of codegree}) and (\ref{lowerbound of edges}), for $1\leq i<j\leq r$, we have 
\[
    \frac{e(G)}{\p}\leq e(H') \leq |\partial_{ij}|\Delta_{ij}\leq n^2\Delta_{ij}.
\]
Thus, we have 
\begin{equation}\label{lower bound of Delta}
    \Delta_{ij}\geq \frac{n^{-2}e(G)}{\p}.
\end{equation}
By (\ref{lower bound of Delta}) , (\ref{extend}) is at least

\[\left(\frac{\Delta_{12}}{4(2\log{n^r})^{r^2}}\right)^{2\ell-1} \left(\frac{\Delta_{1a}}{4(2\log{n^r})^{r^2}}\right)\left(\frac{\Delta_{2a}}{4(2\log{n^r})^{r^2}}\right)=\Delta_{12}^{2\ell-1}\Delta_{1a}\Delta_{2a}\p. \]

  Let $\Ext(C')$ be the collection of all cycles $C_{2\ell+1}^{(r)}$ obtained from $C'$ in this manner described above. Define $\mathcal{F}' =\{\Ext(C')|C\in \mathcal{F}\}$, so $\mathcal{F}'$ is a collection of copies of linear $C_{2\ell+1}^{(r)}$ in $H'$. By the previous discussion and $|\mathcal{F}|\geq \delta_0 k^{2\ell}m^2$, we get
$|\mathcal{F}'|\geq \delta_0 k^{2l}m^2\Delta_{12}^{2\ell-1+\frac{1}{r-2}}\Delta_{1a}\Delta_{2a}\p$.
%
%
%
Now putting all estimates together, we obtain
\begin{align*}
        |\mathcal{F}'| & \geq\delta_0 k^{2l}m^2  \Delta_{12}^{2\ell-1+\frac{1}{r-2}}\Delta_{1a}\Delta_{2a}\p \\
&\geq\delta_0|\partial_{12}|^{2\ell}m^{-2\ell}\Delta_{12}^{2\ell-1+\frac{1}{r-2}}n^{-4}e(G)^2\p \quad\quad\quad\text{ by $k=\tfrac{|\partial_{12}|}{m^{1+1/\ell}}$ and \eqref{lower bound of Delta}}\\
       & \geq\p (|\partial_{12}|\Delta_{12})^{2\ell-1+\frac{1}{r-2}}|\partial_{12}|^{1-\frac{1}{r-2}}n^{-4}e(G)^2m^{-2\ell}\\
       & \geq  \p  e(G)^{2\ell-1+\frac{1}{r-2}} (e(G)^{\frac{1}{r-1}}m^{\frac{r-2}{r-1}})^{1-\frac1{r-2}} n^{-4}e(G)^2m^{-2\ell} \quad\quad \text{ by \eqref{eq:partialedges}}\\
       & \ge \p  e(G)^{2\ell+1+\frac{2}{r-1}} n^{-2\ell-4-\frac{r-3}{r-1}} \quad\quad\quad\quad\quad\quad\text{ by $m\le n$}\\
       & \ge \p  n^{(r-1)(2\ell+1)}\l(\frac{e(G)}{n^r}\r)^{(2\ell+1+\frac{2}{r-1})}\\
       & \ge n^{(r-1)(2\ell+1)-o(1)}\l(\frac{e(G)}{n^r}\r)^{(2\ell+1+\frac{2}{r-1})}.
\end{align*}


For the "as a consequence" part of the theorem, by Proposition~\ref{Copies of F}, setting $e(G)=n^{r-\delta}$, we have 
\[
n^{(r-1)(2\ell+1)-o(1)-\delta(2\ell+1+\frac{2}{r-1})}\leq n^{(r-1)(2\ell+1)-\delta(2\ell+1+s(C^r_{2\ell+1})-\varepsilon)}.
\]
Thus, 
\[
s(C^r_{2\ell+1})\leq\frac{2}{r-1}.
\]
This completes the proofs.
\end{proof}

\subsection{Proof of Theorem~\ref{theorem:main_upperbound}}
We will use induction on $r$ to prove the supersaturation results for $r\ge 4$. The proof splits into two cases according to how large the typical codegrees are: when the typical codegree is small, we find a large collection of $(r-1)$ shadows, use inductive hypothesis on them to obtain a large collection of $C^{r-1}_{2\ell+1}$, and then expand them into a large collection of $C^r_{2\ell+1}$; on the other hand, if the typical codegree is large, then we use the method of ``Greedy Expansion''. We need the following codegree dichotomy lemma.

\begin{lemma}[\cite{nie2024random}, Lemma 3.4]\label{dichotomy}
For $r\ge 2$, let $H$ be an $r$-partite $r$-graph on $n$ vertices, then for any $e(H)/(4n^{r-1})< A\le n$, one of the following two statements is true: 
    \begin{itemize}
        \item [(i)] There exists a subgraph $\hat{H}\subseteq H$ with $e(\hat{H})\ge e(H)/2$ such that every $(r-1)$-shadow of $\hat{H}$ has codegree at least $A$.
        \item[(ii)] There  exist a subgraph $\hat{H}\subseteq H$ with $e(\hat{H})\ge e(H)/(4r\log n)$, an $r$-partition $V(\hat{H})=V_1\cup\dots\cup V_r$, and
        $$
        \frac{e(H)}{4n^{r-1}}<D\le A
        $$
        such that any $(r-1)$-shadow $\sigma$ of $\hat{H}$ in $V_1\times\dots\times V_{r-1}$ satisfies 
        $$
        D/2\le d_{\hat{H}}(\sigma)<D.
        $$
    \end{itemize}
\end{lemma}

For convenience, we restate \Cref{theorem:main_upperbound}.

\begingroup
\def\thethm{\ref{theorem:main_upperbound}}
\begin{thm}
For every $r\geq3$, $\ell\geq 2$, there exists a constant $C$ such that if $G$ is an $n$-vertex $r$-graph such that $e(G)\ge n^{r-1}(\log n)^{C}$, then, as $n\rightarrow\infty$, $G$ contains at least
\[
n^{(r-1)(2\ell+1)-o(1)}\l(\frac{e(G)}{n^{r}}\r)^{\l(2\ell+1+\frac{4\ell-1}{(r-1)(2\ell+1)-3}\r)}
\]
copies of $C^r_{2\ell+1}$. As a consequence,
we have 
\[
s(C^r_{2\ell+1})\le \frac{4\ell-1}{(r-1)(2\ell+1)-3}.
\]
\end{thm}
\endgroup

\begin{proof}[Proof of Theorem~\ref{theorem:main_upperbound}]
     We argue by induction on $r$. The base case $r=3$ follows from Proposition~\ref{theorem:weaker main_upperbound}, set $f(r)=\frac{4\ell-1}{(r-1)(2\ell+1)-3}.$ Assume that for $r-1$,  we have already obtained the desired result. 
     
    Let $G$ be a given $r$-graph with $n$ vertices such that $e(G)\geq n^{r-1}(\log n)^{C}$, where $C$ is a sufficiently large constant.  
     According to a classical result by Erd\H{o}s and Kleitman~\cite{erdos1968coloring}, $G$ contains an $r$-partite subgraph $H$ with $r$-partition $(V_1,V_2,\dots, V_r)$
such that the number of edges in $H$ satisfies $e(H)\geq r!e(G)/r^r$.
Apply Lemma~\ref{dichotomy} on $H$ with \[A=e(G)^{\frac{2\ell+f(r-1)}{2\ell(r-1)-1+f(r-1)}}n^{-\frac{2\ell+1+(r-1)f(r-1)}{2\ell(r-1)-1+f(r-1)}-o(1)}<n^{\frac{2\ell(r-2)-1}{2\ell(r-1)-1+f(r-1)}-o(1)}<n.\] Let $\hat{H}$ be the subgraph of $H$ guaranteed by Lemma~\ref{dichotomy}.

If $(i)$ is true, then there exists a subgraph $\hat{H}\subseteq H$ with $e(\hat{H})\geq e(H)/2\geq r!e(G)/2r^r=\Omega(e(G))$, such that every $(r-1)$-shadow of $\hat{H}$ has codegree at least $A>2(2\ell+1)(r-1)$. 
Now we describe an algorithm to construct copies of $C^r_{2\ell+1}$ in $\hat{H}$; 
it involves specifying $2\ell+1$ edges $e_1, e_2, \dots , e_{2\ell+1}$ which form a copy of $C^r_{2\ell+1}$ in $\hat{H}$, where $e_i=\{v_{i,1}, \dots, v_{i,r-2}, w_{i-1}, w_i\}$, $1\leq i\leq 2\ell+1$, $w_0=w_{2\ell+1}$. 
In this algorithm, we can view an edge as an ordered set of vertices. Given an edge with an order on its vertices $e=(a_1, a_2, \dots, a_r)$ and a vertex $v$, let $e[2,r]=(a_2,\dots, a_r)$ and let $e[2,r]\oplus \{v\}=(a_2,\dots, a_r, v)$.
\begin{itemize}
    \item[(\uppercase\expandafter{\romannumeral1})] Start by choosing an edge with an order $e_1=\{v_{1,1}, \dots, v_{1,r-2}, w_{2\ell+1}, w_1\}$ in $E(\hat{H})$;
    \item [(\uppercase\expandafter{\romannumeral2})] For $1\leq i\leq \ell-1$, we specify $w_{2\ell+1-i}$ and $w_{i+1}$ inductively as following. Let $f_1=e_1$. Next we greedily pick $f_{2\ell+1}=f_1[2,r]\oplus\{w_{2\ell}\}$ such that $w_{2\ell}$ is different from any specified vertices. Then we pick $f_{2}=f_{2\ell+1}[2,r]\oplus\{w_2\}$ such that $w_2$ is different from any other specified vertices. In general, given $f_i$, we choose $f_{2\ell+2-i}=f_i[2,r]\oplus\{w_{2\ell+1-i}\}$ such that the only vertex $w_{2\ell+1-i}$ in $f_{2\ell+2-i}\setminus f_i$ is distinct from all specified vertices. Then we choose $f_{i+1}=f_{2\ell+2-i}[2,r]\oplus\{w_{i+1}\}$ such that the only vertex $w_{i+1}$ in $f_{i+1}\setminus f_{2\ell+2-i}$ is distinct from all specified vertices. Further, given $f_{\ell}$, we choose $f_{\ell+1}=f_{\ell}[2,r]\oplus\{w_{\ell+1}\}$ such that the only vertex $w_{\ell+1}$ in $f_{\ell+1}\setminus f_{\ell}$ is distinct from all specified vertices. Then we let $f_{\ell+2}=f_{\ell+1}$. So far we have $\{w_{i-1},~w_i\}\subset f_i$ for every $2\le i\le 2\ell+1$. (See \Cref{fig:greedy} for an illustration)

\begin{figure}[h]
    \centering
    \includegraphics[width=0.5\linewidth]{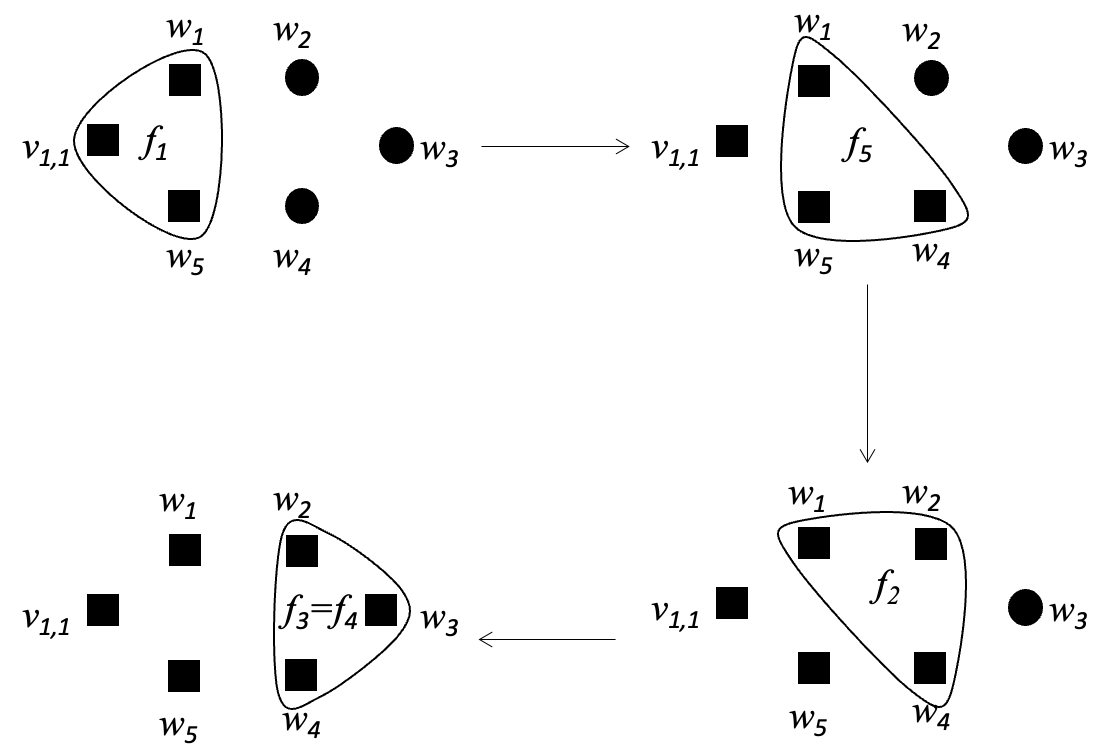}
    \caption{Step (II) of the algorithm for $C^3_5$. The $\blacksquare$s are specified vertices.}
    \label{fig:greedy}
\end{figure}

    \item[(\uppercase\expandafter{\romannumeral3})] For $2\le i\le 2\ell+1$, we give $f_i$ a new order such that $w_{i-1}$ and $w_i$ are the last two vertices. Then we specify $v_{i,j}$, $1\le j\le r-2$, inductively as following. Let $e_{i,0}=f_i$. We choose $e_{i,j}$ such that the only vertex $v_{i,j}$ in $e_{i,j}\setminus e_{i,j-1}$ is distinct from all specified vertices; give vertices of $e_{i,j}$ an order such that $e_{i,j}=e_{i,j-1}[2,r]\oplus\{v_{i,j}\}$. Let $e_{i}=e_{i,r-2}$.
    
\end{itemize}

Note that $e_1, e_2, \dots, e_{2\ell+1}$ form a copy of $C^r_{2\ell+1}$.In the above algorithm, in the first step, we have $e(\hat{H})$ choices for selecting $e_1$. 
In the subsequent steps, each time we have at least $A-(2\ell+1)(r-1)$ choices for selecting a new vertex, since every $(r-1)$-shadow of $\hat{H}$ has codegree at least $A$. 
Recall that $A> 2(2\ell+1)(r-1)$.
Thus, we can find at least $\Omega(e(G)A^{2\ell(r-1)-1})$ labeled copies of $C^r_{2\ell+1}$ in $\hat{H}$. Note that
\begin{align*}
    \begin{split}
        e(G)A^{2\ell(r-1)-1}&\geq e(G)\cdot e(G)^{\frac{(2\ell+f(r-1))(2\ell(r-1)-1)}{2\ell(r-1)-1+f(r-1)}}n^{-\frac{(2\ell+1+(r-1)f(r-1))(2\ell(r-1)-1)}{2\ell(r-1)-1+f(r-1)}-o(1)}\\
        &= e(G)\cdot e(G)^{\frac{(2\ell[(r-2)(2\ell+1)-3]+4\ell-1)(2\ell(r-1)-1)}{(2\ell(r-1)-1)[(r-2)(2\ell+1)-3]+4\ell-1}}n^{-(2\ell+1)-\frac{(r-1)f(r-1)(2\ell(r-1)-1)-(2\ell+1)f(r-1)}{2\ell(r-1)-1+f(r-1)}-o(1)}\\
        &= e(G)^{2\ell+1+\frac{(4\ell-1)(2\ell(r-1)-1)-2\ell(4\ell-1)}{(2\ell(r-1)-1)[(r-2)(2\ell+1)-3]+4\ell-1}}n^{-(2\ell+1)-\frac{(r-1)(4\ell-1)(2\ell(r-1)-1)-(2\ell+1)(4\ell-1)}{(2\ell(r-1)-1)[(r-2)(2\ell+1)-3]+4\ell-1}-o(1)}\\
        &= e(G)^{2\ell+1+\frac{(4\ell-1)(2\ell(r-2)-1)}{(2\ell(r-2)-1)[(r-1)(2\ell+1)-3]}}n^{-(2\ell+1)-\frac{r(4\ell-1)(2\ell(r-2)-1)}{(2\ell(r-2)-1)[(r-1)(2\ell+1)-3]}-o(1)}\\
        &=e(G)^{2\ell+1+\frac{4\ell-1}{(r-1)(2\ell+1)-3}}n^{-(2\ell+1)-\frac{r(4\ell-1)}{(r-1)(2\ell+1)-3}-o(1)}\\
        &=n^{(r-1)(2\ell+1)-o(1)}\l(\frac{e(G)}{n^{r}}\r)^{\l(2\ell+1+\frac{4\ell-1}{(r-1)(2\ell+1)-3}\r)}.
    \end{split}
\end{align*}
Thus we are done in this case.

If $(ii)$ is true, then there exist a subgraph $\hat{H}\subseteq H$ with $e(\hat{H})\ge e(H)/(4r\log n) \geq C_2e(G)/\log n$, where $C_2$ is a constant, an $r$-partition $V(\hat{H})=V_1\cup\dots\cup V_r$, and $\frac{e(H)}{4n^{r-1}}<D\le A$
such that any $(r-1)$-shadow $\sigma$ of $\hat{H}$ in $V_1\times\dots\times V_{r-1}$ satisfies $D/2\le d_{\hat{H}}(\sigma)<D.$ 
Consider the $(r-1)$-graph $H'$ on $V_1\times\dots\times V_{r-1}$, whose edge set contains all the $(r-1)$-shadows of $\hat{H}$ in $V_1\times\dots\times V_{r-1}$. Note that we have $e(H')\cdot D\ge e(\hat{H})$. Thus, $e(H')\ge \frac{e(\hat{H})}{D}=C_2e(G)(\frac{(\log n)^{-1}}{D})\geq C_2e(G)(\frac{(\log n)^{-1}}{n})\geq C_2n^{r-2}(\log n)^{C-1}$.
 Then by inductive hypothesis, we have the number of $C^{r-1}_{2\ell+1}$ copies in $H'$ is at least $n^{(r-2)(2\ell+1)-o(1)}\left(\frac{e(G)n^{-o(1)}}{n^{r-1}D}\right)^{(2\ell+1)+f(r-1)}$. 
Moreover, the number of $C^r_{2\ell+1}$ copies in $\hat{H}$ is at least 
\begin{equation}\label{copies of C}
n^{(r-2)(2\ell+1)-o(1)}\left(\frac{e(G)n^{-o(1)}}{n^{r-1}D}\right)^{(2\ell+1)+f(r-1)}\cdot D^{2\ell+1},
\end{equation}
since each $(r-1)$-shadow of $\hat{H}$ can be extended into an $r$-edge of $\hat{H}$ in at least $D/2-2\ell$ ways. 
In addition, ~\eqref{copies of C} is equal to 
\[
n^{-(2\ell+1)-(r-1)f(r-1)-o(1)}e(G)^{2\ell+1+f(r-1)}D^{-f(r-1)}\geq n^{-(2\ell+1)-(r-1)f(r-1)-o(1)}e(G)^{2\ell+1+f(r-1)}A^{-f(r-1)}.
\]
Note that
\begin{align*}
    \begin{split}
        RHS&=n^{(r-1)(2\ell+1)-o(1)}\l(\frac{e(G)}{n^{r}}\r)^{2\ell+1}n^{\frac{f(r-1)[(2\ell+1)-2\ell(r-1)^2+(r-1)]}{2\ell(r-1)-1+f(r-1)}}e(G)^{\frac{f(r-1)[2\ell(r-1)-1-2\ell]}{2\ell(r-1)-1+f(r-1)}}\\
        &=n^{(r-1)(2\ell+1)-o(1)}\l(\frac{e(G)}{n^{r}}\r)^{2\ell+1}n^{\frac{-r(4\ell-1)(2\ell(r-2)-1)}{(2\ell(r-2)-1)[(r-1)(2\ell+1)-3]}}e(G)^{\frac{(4\ell-1)(2\ell(r-2)-1)}{(2\ell(r-2)-1)[(r-1)(2\ell+1)-3]}}\\
        &=n^{(r-1)(2\ell+1)-o(1)}\l(\frac{e(G)}{n^{r}}\r)^{2\ell+1}n^{\frac{-r(4\ell-1)}{(r-1)(2\ell+1)-3}}e(G)^{\frac{4\ell-1}{(r-1)(2\ell+1)-3}}\\
        &\geq n^{(r-1)(2\ell+1)-o(1)}\l(\frac{e(G)}{n^{r}}\r)^{(2\ell+1+\frac{4\ell-1}{(r-1)(2\ell+1)-3})}.
    \end{split}
\end{align*}
Thus we are also done in this case.


For the "as a consequence" part of the theorem, by Proposition~\ref{Copies of F}, setting $e(G)=n^{r-\delta}$, we have 
\[
n^{(r-1)(2\ell+1)-o(1)-\delta(2\ell+1+\frac{4\ell-1}{(r-1)(2\ell+1)-3})}\leq n^{(r-1)(2\ell+1)-\delta(2\ell+1+s(C^r_{2\ell+1})-\varepsilon)}.
\]
Thus, 
\[
s(C^r_{2\ell+1})\leq\frac{4\ell-1}{(r-1)(2\ell+1)-3}.
\]
This completes the proofs.
\end{proof}

\section{Conclusion}
\begin{itemize}
    \item We believe the conditional lower bound in \Cref{theorem:main_lowerbound} gives the correct value of $s(C^r_{2\ell+1})$. The current approach in this paper falls short in proving a matching upper bound. New ideas will be necessary to improve the upper bound in \Cref{theorem:main_upperbound}.
    \begin{conjecture}
    For $r\ge 3$ and $\ell\ge 2$, $s(C^r_{2\ell+1})=\frac{1}{(r-1)\ell-1}$.
    \end{conjecture}
    
    \item Note that in the proof of the Proposition~\ref{theorem:weaker main_upperbound}, we used the supersaturation result for even cycles from~\cite{morris2016number}, which is actually a balanced supersaturation result. 
Indeed, using the full power of this result, we can also derive a hypergraph version of the balanced supersaturation result for linear odd cycles. 
By applying the container method, we can obtain the random Tur{\'a}n result for $3$-uniform linear odd cycles as follows:

\begin{theorem}\label{thm:random turan}For every $\ell\geq 2$  a.a.s. the following holds:
$$\ex\left(G_{n,p}^{3}, C_{2\ell+1}^{3}\right)\leq \begin{cases}
    p^{\frac{1}{2\ell+1}}n^{1+\frac{3}{2\ell+1}+o(1)}, & \mbox{if } {n^{-2+\frac{1}{2\ell}}\leq p\ll n^{-2+\frac{2}{2\ell-1}}} \\
    p^{\frac{1}{2}}n^{2+o(1)}, & \mbox{if } p\gg n^{-2+\frac{2}{2\ell-1}}.
\end{cases}
$$
\end{theorem}
However, the proof of this result (using the container method) is involved and by now standard. Therefore, we choose to omit the complex proof of a suboptimal result.
    \item One can also consider a stronger version of \Cref{conjecture:supersat} by replacing the $n^{-o(1)}$ factor by a constant. Such a stronger version is (implicitly) known to be true~\cite{balogh2019number} for $C^3_3$. But the current results for $C^r_{2\ell}$ and the $r$-expansion of $K^k_{k+1}$ are insufficient for this stronger version. It would be interesting to prove this stronger version of \Cref{conjecture:supersat} for any hypergraph other than $C^3_3$.
\end{itemize}

\bibliographystyle{plain}
\bibliography{ref}

\end{document}